\documentclass[a4paper,11pt]{article} 

\usepackage{a4wide}
\usepackage{pdflscape}
\usepackage{makeidx}
\usepackage{version}
\usepackage{amsmath,amssymb,amsfonts,amsthm,amscd}
\usepackage{color, pgf, latexsym}
\usepackage{float} 
\usepackage{algpseudocode}
\usepackage{algorithm}
\usepackage{tikz}
\usepackage{listings,setspace} 

\input{xy}
\xyoption{all}
\input xy \xyoption{matrix} \xyoption{arrow}
\def\edge{\ar@{-}}

\newtheorem{theorem}{Theorem}[section]
\newtheorem{corollary}[theorem]{Corollary}
\newtheorem{lemma}[theorem]{Lemma}
\newtheorem{proposition}[theorem]{Proposition}

\theoremstyle{definition}
\newtheorem{definition}[theorem]{Definition}

\newtheorem{remark}[theorem]{Remark}

\title{Factorisation of cross-symmetric, totally nonnegative matrices and an amazing matrix}
\author{T H Lenagan and A P Neate\footnote{The work of the second author was supported by a University of Edinburgh School of Mathematics Summer Vacation Scholarship}}
\date{}

\begin{document}
\maketitle
\begin{abstract} We establish a factorisation theorem for invertible, cross-symmetric, totally nonnegative matrices, and illustrate the theory by verifying that certain cases of Holte's Amazing Matrix are totally nonnegative. 
\end{abstract}

\vskip .5cm
\noindent
{\bf 2020 Mathematics subject classification:} 15B48
\vskip .5cm
\noindent
\textbf{Keywords.} Cross-symmetric matrix, centro-symmetric matrix,  totally nonnegative matrix

\section{Introduction} 
In \cite{holte}, Holte introduces an ``Amazing Matrix" that arises in connection with carries that occur when adding integers. For each base $b$, he constructs an $n\times n$ matrix $P$ whose $(i,j)$-entry is the probablility that, when adding $n$ random numbers written in base $b$, the next carry will be $j$ given that the previous carry was $i$. Holte discovered many interesting facts about the amazing matrix, including the fact that $P$ is cross-symmetric (centro-symmetic in some sources); that is, it is radially symmetric about its centre.

In a later article \cite{df}, Diaconis and Fulman uncovered a connection between the amazing matrix and card shuffling. In their study of properties of the amazing matrix, they showed that when the base $b=2$ the matrix is totally nonnegative; that is, all minors are nonnegative. They also showed that in arbitrary base $b$ all of the $2\times 2$ minors are nonnegative, and they made the conjecture that the amazing matrix is totally nonnegative. As far as we are aware, this conjecture is still open. In \cite{mcmillan}, McMillan made some progress towards verifying that the conjecture holds. 

There are efficient tests to determine whether or not a matrix is totally nonegative, for example, Neville elimination \cite[Chapter 2]{fj} or Cauchon's deletion algorithm \cite{gll}. 
However, these processes immediately destroy the desirable property of cross-symmetry. Our interest in the amazing matrix conjecture led us to construct the algorithm presented here which preserves cross-symmetry and checks for total nonnegativity. It leads to the factorisation theory that we present here: each invertible, cross-symmetric, totally nonnegative matrix can be factorised into a product of cross-symmetric ``atoms'' which are the cross-symmetric analogues of the elementary matrices that arise in the Neville factorisation theory. 

Although we are not able to verify the amazing matrix conjecture in general, we can use the algorithm to test for specific values of $n$. As an example, the algorithm has been run on amazing matrices of size less than or equal to six and shows that in arbitrary base $b$  all such matrices are totally nonnegative.


\section{Definitions and basic results} 

A matrix with real number entries is said to be {\em totally nonnegative} if each of its minors is nonnegative. Totally nonnegative matrices arise in many different settings, for example, oscillations in mechanical
systems, stochastic processes and approximation theory, P\'olya frequency sequences,
representation theory, planar networks, ... . Two recent books which serve as useful references for properties of totally nonnegative matrices are \cite{fj} and \cite{pinkus}. (The reader should be aware that, in \cite{pinkus}, Pinkus uses the term {\em totally positive} where we use the term totally nonnegative.) \\

If $M$ is a matrix and $I$ and $J$ are sets of row and column indices then we will denote by $M(I,J)$ the submatrix of $M$ formed by using rows $I$ and $J$. Also, if $I$ and $J$ have the same size then we set the minor $[I\mid J]_M$ to be the determinant of $M(I,J)$. If it is obvious which matrix is being discussed, then we may drop the subscript $M$ and simply write $[I\mid J]$.\\

Let $w_0$ be the longest element of the permutation group $S_n$; so that $w_0(i)=n+1-i$, for each $1\leq i\leq n$. An $n\times n$ matrix $M=(m_{ij})$ is {\em cross-symmetric} (or {\em centro-symmetric} in some sources) if $m_{w_0(i),w_0(j)}=m_{ij}$, for all $1\leq i,j\leq n$. Basic properties of cross-symmetric matrices are discussed in \cite{weaver}.
\\

Define $\tau(M)$ to be the matrix whose $(i,j)$ entry is given by $m_{w_0(i),w_0(j)}$. Note that $M$ is cross-symmetric if and only if $\tau(M)=M$. 

Let $J$ be the $n\times n$ matrix whose entries in positions $(i,n+1-i)$ are all equal to $1$, while all other entries are $0$; that is, $J$ has $1$ in each position on the anti-diagonal and is zero elsewhere. It is easy to check that $\tau(M) = JMJ$. Thus, $M$ is cross-symmetric if and only if $M=JMJ$. It follows easily that 
$\tau(AB)=\tau(A)\tau(B)$ and that $AB$ is cross-symmetric whenever each of $A$ and $B$ is cross-symmetric. 

First, we want to consider the minors of the matrix $\tau(A)$ in terms of the minors of $A$. 

\begin{lemma}
Let $I$ and $J$ be index sets of the same size. Then 
\[
[I\mid J]_{\tau(A)}=[w_0(I)\mid w_0(J)]_A.
\]
\end{lemma}

\begin{proof} 
When $|I|=|J|=1$ this is the definition of $\tau(A)$. The proof is then by induction on $|I|$ using induction and Laplace expansions. 
\end{proof}

\begin{corollary}\label{corollay-tau-tnn}
The matrix $A$ is totally nonnegative if and only if $\tau(A)$ is totally nonnegative.
\end{corollary} 

\begin{lemma}
Let $A$ be a totally nonnegative matrix. Suppose that there is a pair $(s,t)$ such that $a_{sj}=0$ for $j\leq t$, but $a_{s+1,t}\neq 0$. Then $a_{sj}=0$ for all $j.$ 
\end{lemma}

\begin{proof}
Let $j>t$. It  is enough to show that $a_{sj}=0$.  As $A$ is totally nonnegative, 
\[
0\leq [s,s+1\mid t,j] = a_{st}a_{s+1,j}-a_{s+1,t}a_{sj}=0-a_{s+1,t}a_{sj}\leq 0,
\]
so that $a_{s+1,t}a_{sj}=0$
It follows that $a_{sj}=0$, as $a_{s+1,t}>0$.
\end{proof}

\begin{corollary}\label{corollary-no-zeroes} 
Let $A$ be an invertible, totally nonnegative matrix. Then there is no pair $(s,t)$ such that $a_{sj}=0$ for $j\leq t$, but $a_{s+1,t}\neq 0$.
\end{corollary} 

\begin{proof}
If $A$ had such a pair $(s,t)$ then the previous lemma would imply that row $s$ of $A$ is zero, which is impossible for an invertible matrix.
\end{proof} 

The following result is a special case of \cite[Theorem 1.13]{pinkus}
\begin{lemma}\label{lemma-diagonals-positive} 
The diagonal elements of an invertible, totally nonnegative matrix are all greater than zero.
\end{lemma}

\section{Cross-symmetric elimination} 
First, we give an informal description of the elimination process that we will use. We start with an invertible, cross-symmetric, totally nonnegative matrix. The aim is to use a version of the Neville elimination procedure to produce a final matrix that is a diagonal matrix. We proceed as with Neville elimination: if we are clearing the lower entries in a given column and want to perform a row operation to replace the last nonzero entry in a column by zero, then we perform a row operation by subtracting a suitable multiple of the row immediately above this last position. Note that, if the matrix is totally nonnegative, the entry immediately above this last position will be nonzero: this is guaranteed by Corollary~\ref{corollary-no-zeroes}. Suppose that this operation involves rows $s$ and $s+1$, then in order to preserve cross-symmetry, at the same time we need to subtract a suitable multiple of row $w_0(s)$ from row $w_0(s+1)$. We need to show that total nonnegativity is preserved by such a move. A complication arises when $n$ is even and $s=n/2$ (which happens precisely when $w_0(s)=s+1$),  so we treat this case separately in the next results where we show what happens under such moves. \\

We denote by $E(i,j)$ the matrix that is zero in all positions except the $(i,j)$ position, where the entry is $1$. An {\em elementary cross-symmetric matrix} is a matrix of the form 
\[
F:=I-cE(s+1,s)-cE(w_0(s+1),w_0(s))
\]
for some  $c>0$. 

The reduction step for our algorithm will consist of pre-multiplying an invertible, cross-symmetric matrix $A$ by a suitable elementary cross-symmetric matrix $F$. It is obvious that such an $F$ is invertible and cross-symmetric and so $FA$ is also invertible and cross-symmetric. It will remain to show that $FA$ is totally nonnegative  provided that $A$ is totally nonegative. When $n\neq 2s$ this is easy to do: $w_0(s+1)=(n+1)-(s+1)=n-s\neq s$ and so $E(s+1,s)E(w_0(s+1),w_0(s))=0$. It follows that $F=\left(I-cE(s+1,s)\right)\left(I-cE(w_0(s+1),w_0(s))\right)$. If we write these two factors as $F_1$ and $F_2$, we see that $FA=F_1F_2A$, and we can control total nonegativity by proceeding in two stages: first premultiply by $F_2$ and then multiply by $F_1$ and use \cite[Proposition 2.6]{gl} twice; precise details are given below. 

\begin{remark} 
\cite[Proposition 2.6]{gl} is essentially due to Whitney \cite{whitney}, see \cite[Theorem 2.2.1]{fj}.
\end{remark}

In order to deal with the first of these products, we need the following preparatory lemma.

\begin{lemma} \label{lemma-half-tnn}
Suppose that $n\neq 2s$ and that $A$ is a cross-symmetric,  totally nonnegative $n\times n$ matrix with $a_{st}\neq 0$, while $a_{ij}=0$ whenever $i\geq s$ and $j<t$. Suppose that $a_{s+1,t}\neq 0$ while $a_{s+w,t}=0$ for all $w>1$. Set 
\[
F:=I-a_{w_0(s+1),w_0(t)}a_{w_0(s)w_0(t)}^{-1}E(w_0(s+1),w_0(s)).
\]
 Then\\ (i) $B=(b_{ij}):=FA$ is totally nonnegative, and\\ (ii)  $b_{st}=a_{st}\neq 0$, while $b_{ij}=0$ for $i\geq s$ and $j<t$. Also, $b_{s+1,t}=a_{s+1,t}\neq 0$ while $b_{s+w,t}=0$ for all $w>1$.
\end{lemma}

\begin{proof} (i) As $\tau(F)= I-a_{s+1,t}a_{st}^{-1}E(s+1,s)$ the pair of matrices $\tau(F)$ and $A$ satisfy the conditions necessary to apply  \cite[Lemma 2.5, Proposition 2.6]{gl} and from that to conclude that $\tau(FA)=\tau(F)\tau(A)$ is totally nonnegative. 
It follows that $\tau(B)=\tau(FA)=\tau(F)\tau(A)=\tau(F)A$ is totally nonnegative. As a result, $B$ is totally nonnegative, by Corollary~\ref{corollay-tau-tnn}.\\

(ii) Note that $n\neq 2s$ implies that $w_0(s+1)= n+1-(s+1)=n-s\neq s$. The matrices $FA$ and $A$ are identical in all rows except row $w_0(s+1)$. Hence, if  $w_0(s+1)<s$ then rows $s+i$ with $i\geq 0$ are the same in $FA$ and $A$ and so the claims in statement (ii) of the lemma follow.

Now, suppose that $w_0(s+1)>s$. The rows other than row $w_0(s+1)$ are the same for $B$ and $A$, so we only have to check the claims in statement (ii) for the relevant entries of $B$ on row $w_0(s+1)$, where, to obtain the value of $b_{w_0(s+1),j}$ we subtract from $a_{w_0(s+1),j}$ a multiple of $a_{w_0(s),j}$. 
However, $w_0(s)>w_0(s+1)>s+1$, and so $w_0(s)\geq s+2$. It follows that $a_{w_0(s),j}=0$ for $j\leq t$, and again we conclude that the claims in statement (ii) of the lemma follow. 
\end{proof}

In order to deal with the case where $n=2s$, we need another preparatory lemma.

\begin{lemma}\label{lemma-n=2s}
Suppose that $n=2s$ and that $A$ is an invertible,  cross-symmetric, totally nonnegative $n\times n$ matrix with $a_{st}\neq 0$, while $a_{ij}=0$ whenever $i\geq s$ and $j<t$. Suppose that $a_{s+1,t}\neq 0$ while $a_{s+w,t}=0$ for all $w>1$. Then $t\leq s$, $t<w_0(t)$ and $a_{s+1,t}<a_{st}$

\end{lemma}

\begin{proof}
In this case, $w_0(s)=n+1-s=s+1$ and $w_0(s+1)=s$. 
Suppose that  $t>s$ then the hypotheses show that $a_{ss}=0$, which contradicts Lemma~\ref{lemma-diagonals-positive}. It follows that $t\leq s$, and that  $t\leq s<s+1=w_0(s)\leq w_0(t)$.

As $A$ is cross-symmetric, $a_{s+1,w_0(t)}=a_{w_0(s),w_0(t)}=a_{st}$ and $a_{s,w_0(t)}=a_{w_0(s+1)w_0(t)}=a_{s+1,t}$. Hence, as $A$ is totally nonnegative, 
\[
0\leq[s,s+1\mid t,w_0(t)] = a_{st}a_{s+1,w_0(t)}-a_{s,w_0(t)}a_{s+1,t} = a_{st}^2-a_{s+1,t}^2\,,
\]
so that $a_{s+1,t}\leq a_{st}$. 

Suppose that $a_{s+1,t}=a_{st}$. As $A$ is cross-symmetric,  $a_{w_0(s+1),w_0(t)}=a_{w_0(s),w_0(t)}$; that is, 
$a_{s,w_0(t)}=a_{s+1,w_0(t)}$. 
 First, suppose that $w_0(t)=t+1$. Then rows $s$ and $s+1$ are equal, since they both take the form $0,\dots,0,a_{st},a_{s,t},0,\dots,0$. This cannot happen, as $A$ is invertible. If $w_0(t)>t+1$, then consider any $j$ such that $t<j<w_0(t)$. We know that $0\leq[s,s+1\mid t,j]=a_{st}a_{s+1,j}-a_{sj}a_{s+1,t}=
a_{st}a_{s+1,j}-a_{sj}a_{s,t}$. As $a_{st}>0$, this implies that $a_{s+1,j}\geq a_{sj}$. A  similar argument with 
$[s,s+1\mid j,w_0(t)]$ shows that $a_{sj}\geq a_{s+1,j}$, and so $a_{sj}=a_{s+1,j}$. Again, we see that rows $s$ and $s+1$ are equal, a contradiction; so that $a_{s+1,t}<a_{st}$, as required. 
\end{proof}

\begin{proposition} 
Suppose that  $A=(a_{ij})$ is an $n\times n$ invertible, cross-symmetric, totally nonnegative matrix with $a_{st}\neq 0$ and $a_{ij}=0$ whenever $i\geq s$ and $j<t$. Suppose also that $a_{s+1,t}\neq 0$ while $a_{s+w,t}=0$ for all $w>1$. Set 
\[
F:=I-a_{s+1,t}a_{st}^{-1}E(s+1,s)-a_{w_0(s+1),w_0(t)}a_{w_0(s)w_0(t)}^{-1}E(w_0(s+1),w_0(s))
\]
and set $B=(b_{ij}):=FA$. Then $B$ is an invertible, cross-symmetric, totally nonnegative matrix.

Further, $b_{st}\neq 0$, and $b_{ij}=0$ for $i\geq s$ and $j<t$, while $b_{s+w,t}=0$ for all $w\geq 1$.
\end{proposition}

\begin{proof} As $F$ is invertible and cross-symmetric, we know that $B=FA$ is invertible and cross-symmetric, so we need to show that $B$ is totally nonnegative. 
\\

First, suppose that $n\neq 2s$ so that 
\[
F=\left(I-a_{s+1,t}a_{st}^{-1}E(s+1,s)\right)
\left(I-a_{w_0(s+1),w_0(t)}a_{w_0(s)w_0(t)}^{-1}E(w_0(s+1),w_0(s))\right)
\]

Lemma~\ref{lemma-half-tnn} shows that $B':=(I-a_{w_0(s+1),w_0(t)}a_{w_0(s)w_0(t)}^{-1}E(w_0(s+1),w_0(s))A$ is totally nonnegative, and the entries of $B'$ satisfy the same conditions as those for $A$ that are given in the hypotheses of the proposition, and, in particular, $b'_{st}=a_{st}$ and $b'_{s+1,t}=a_{s+1,t}$. 

Hence, $FA=(I-a_{s+1,t}a_{st}^{-1}E(s+1,s))B'$  is totally nonnegative by 
\cite[Lemma 2.5,Proposition 2.6]{gl}. Also, $b_{st}=b'_{st}=a_{st}\neq 0$ while $b_{s+1,t}=a_{s+1,t}- a_{s+1,t}a_{st}^{-1}a_{st}=0$.\\

Now, assume that $n= 2s$. \\

The matrices $A$ and $B$ differ only in rows $s,s+1$, so we need to know the entries in these rows of $B$: 
\[
b_{sj} = a_{sj}-a_{s,w_0(t)}a_{s+1,w_0(t)}^{-1}a_{s+1,j}
\quad{\rm and}\quad 
b_{s+1,j}=a_{s+1,j}-a_{s+1,t}a_{st}^{-1}a_{sj}.
\]

Let $I$ and $J$ be row and column index sets of the same size. We need to show that $[I\mid J]_B\geq 0$. In order to do this, we consider four cases depending on which of $s,s+1$ is/are in $I$. 
\\

First, suppose that $s,s+1\not\in I$.  Then $B(I,J)=A(I,J)$ and so $[I\mid J]_B=[I\mid J]_A\geq 0$, as $A$ is totally nonnegative.\\

Secondly, suppose  that $s\not\in I$ but $s+1\in I$. Set $C:=(I-a_{s+1,t}a_{st}^{-1}E(s+1,s))A$. Then $B(I,J)=C(I,J)$. However, $C(I,J)$ is totally nonnegative, by \cite[Lemma 2.5, Proposition 2.6]{gl}; and so $[I\mid J]_B=[I\mid J]_C\geq 0$.\\

Next, suppose that $s\in I$, but $s+1\not\in I$. In this case, set 
\[
C:= (I-a_{w_0(s+1),w_0(t)}a_{w_0(s)w_0(t)}^{-1}E(w_0(s+1),w_0(s))A.
\] Then $C$ is totally nonnegative, by 
Lemma~\ref{lemma-half-tnn}. As $B(I,J)=C(I,J)$ we conclude that $[I\mid J]_B\geq 0$.\\

Finally, suppose that $s,s+1\in I$. \\

Let $A_1$ be the matrix obtained by replacing row $s$ of $A$ by row $s+1$ of A; let $A_2$ be the matrix obtained by replacing row $s+1$ of $A$ by row $s$ of A; and let $A_3$ be the matrix obtained from $A$ by swopping rows $s$ and $s+1$. 

Then, using the fact that the determinant is linear in each row (in particular, in rows $s$ and $s+1$), we obtain 
\begin{align*}
[I\mid J]_B= [I\mid J]_{FA} &=[I\mid J]_A - a_{st}a_{s+1,t}^{-1}[I\mid J]_{A_1}
\\
&~~~ - a_{s+1,t}a_{st}^{-1}[I\mid J]_{A_2}
+ a_{st}a_{s+1,t}^{-1}a_{s+1,t}a_{st}^{-1}[I\mid J]_{A_3}
\end{align*}

Now, $[I\mid J]_{A_1}=[I\mid J]_{A_2}=0$, as 
$A_1$ and $A_2$ have repeated rows that are included in the row set $I$. Also, $[I\mid J]_{A_3}=-[I\mid J]_A$ as $A_3$ is obtained from $A$ by swopping two adjacent rows.

Hence, 
\begin{align*}
[I\mid J]_B
&=
[I\mid J]_A-a_{s,w_0(t)}a_{s+1,w_0(t)}^{-1}a_{s+1,t}a_{st}^{-1}[I\mid J]_A\\
&=a_{s+1,w_0(t)}^{-1}a_{st}^{-1}\left( a_{st}a_{s+1,w_0(t)} - a_{s,w_0(t)}a_{s+1,t} \right)[I\mid J]_A
\end{align*}
As $A$ is totally nonnegative, $[I\mid J]_A\geq 0$. Also, $a_{s+1,w_0(t)}=a_{st}>0$ from the hypotheses and the  fact that $A$ is cross-symmetric. 
Recalling that $t<w_0(t)$, we obtain 
\[
a_{st}a_{s+1,w_0(t)} - a_{s,w_0(t)}a_{s+1,t} = [s,s+1\mid t,w_0(t)]_A\geq 0
\]
 and it follows that 
$[I\mid J]_B\geq 0$. \\

Finally, note that 
\begin{align*}
b_{st}&=a_{st} - a_{s,w_0(t)}a_{s+1,w_0(t)}^{-1}a_{s+1,t}=a_{st} - a_{w_0(s+1),w_0(t)}a_{w_0(s),w_0(t)}^{-1}a_{s+1,t}\\
&= a_{st} - a_{s+1,t}a_{s,t}^{-1}a_{s+1,t}=\left(a_{st}^2-a_{s+1,t}^2\right)a_{st}^{-1} >0, 
\end{align*}
where the inequality occurs because we know that $a_{st}>a_{s+1,t}$, by Lemma~\ref{lemma-n=2s}. 
The other claims for the elements of $B$ follow easily. 
\end{proof}

\begin{theorem} Let $A$ be an $n\times n$ invertible, cross-symmetric, totally nonnegative matrix. Then there 
is a sequence of elementary cross-symmetric matrices $F_1,\dots,F_d$ such that 
$F_d\dots F_1A=D$ where $D=(d_{ij})$ is a cross-symmetric, diagonal matrix with $d_{ii}>0$ for $1\leq i\leq n$.
\end{theorem}

\begin{proof} 
We put the positions of the matrix below the diagonal in the order 
\[
(n,1), \dots, (2,1), (n,2),\dots,(3,2), \dots,(n,n-1);
\]
that is, we move from the bottom of the first column upwards, then up the second column, etc. Use the previous result to set each of these elements to zero. As we do this, the cross-symmetry property ensures that the above diagonal elements also become zero. At the end, we reach a diagonal matrix which is totally nonnegative and, in fact, must have strictly positive entries, as it is invertible. 
\end{proof}

In the setting of this theorem, we see that $A=F_1^{-1}\dots F_d^{-1} D$. We now identify the form of these inverses. 

\begin{definition} 
A {\em (cross-symmetric, totally nonnegative) atom} is a matrix $F$ that takes one of the following two forms. Either $n\neq 2s$ and 
$F:=I+cE(s+1,s)+cE(w_0(s+1),w_0(s))$ for some $c>0$, or, $n=2s$ and $F=(f_{ij})$ where $f_{ss}=f_{s+1,s+1}=1/(1-c^2), f_{s,s+1}=f_{s+1,s} =c/(1-c^2)$ with $0<c<1$, while the other diagonal elements are all equal to $1$ and the other off-diagonal elements are all equal to $0$.  (Note that these matrices $F$ are indeed cross-symmetric and totally nonnegative.) 
\end{definition} 

\begin{lemma} 
(i) Suppose that $n\neq 2s$ and that  $F=I-cE(s+1,s)-cE(w_0(s+1),w_0(s))$ for some $c>0$. Then
$F^{-1} =I+cE(s+1,s)+cE(w_0(s+1),w_0(s))$ is a cross-symmetric, totally nonnegative atom of the first kind defined above.

(ii) Suppose that $n= 2s$ and that  $F=I-cE(s+1,s)-cE(w_0(s+1),w_0(s))$ for some $0<c<1$. Then $F^{-1}$ is a cross-symmetric, totally nonnegative atom of the second type defined above. 
\end{lemma} 

\begin{proof} (i) This is just a matter of checking that 
\[
(I-cE(s+1,s)-cE(w_0(s+1),w_0(s)))(I+cE(s+1,s)+cE(w_0(s+1),w_0(s)))=I.
\]
(ii) The $2\times 2$ submatrix of $F$ formed from rows $s,s+1$ and columns $s,s+1$ has the form 
$\left(\begin{smallmatrix}
1&-c\\
-c&1
\end{smallmatrix}\right)$. The inverse of this matrix is
\[
\begin{pmatrix}
\frac{1}{1-c^2}&\frac{c}{1-c^2}\\
\frac{c}{1-c^2}&\frac{1}{1-c^2},
\end{pmatrix}.
\]
Hence,  $F^{-1}$ is easily seen to be a cross-symmetric totally nonnegative atom of the second kind defined above. 
\end{proof}

\begin{corollary} \label{corollary-factorisation}
Let $A$ be an $n\times n$ invertible, cross-symmetric, totally nonnegative matrix. Then there is a factorisation of $A$ as a product of cross-symmetric totally, nonnegative atoms and a cross-symmetric diagonal matrix with strictly positive diagonal entries. 

Conversely, any product $A_1\dots A_t D$ of cross-symmetric, totally nonnegative atoms $A_i$ and a cross-symmetric  diagonal matrix $D$ with strictly positive diagonal entries is an invertible,  cross-symmetric, totally nonnegative matrix.
\end{corollary} 

\begin{proof}
In the setting of the theorem above, we have $A=F_1^{-1}\dots F_d^{-1} D$, and the inverses that occur are cross-symmetric totally nonnegative atoms by the previous lemma. The converse statement is immediate from earlier comments. \end{proof} 


\section{The algorithm}

The results of the previous section provide us with an algorithmic way of testing when an invertible, cross-symmetric $n\times n$ matrix $A$ is totally nonnegative. Informally, the algorithm proceeds as follows. 

We list the  positions  of a matrix below the diagonal in the order $(n,1), (n-1,1), \dots, (2,1),\\(n,2),\dots, (3,2),\dots, (n,n-1)$; that is, we move from the bottom of the first column upwards, then up the second column, etc.

If all of these entries in $A$ are zero then $A$ is a diagonal matrix, so check if all $a_{ii}>0$, in which case the original $A$ is totally nonnegative, otherwise the original $A$ is not totally nonnegative. 

If there are nonzero elements, then suppose that the first in the list is $a_{s+1,t}$. If $a_{st}\leq 0$ then $A$ is not totally nonnegative. If $a_{st}>0$ then replace $A$ by $FA$ where 
\[
F:=I-a_{s+1,t}a_{st}^{-1}E(s+1,s)1-a_{w_0(s+1),w_0(t)}a_{w_0(s)w_0(t)}^{-1}E(w_0(s+1),w_0(s))
\]
and repeat the above instructions with the new $A$. 

If we keep track of the $F$ that arise then calculate each $F^{-1}$, we recover the factorisation given by Corollary~\ref{corollary-factorisation}.


\section{Applications to Holte's Amazing Matrix}

There is a very close connection between totally nonegative matrices and weighted path matrices of directed planar networks.  In fact, theorems of Lindstr\"om \cite{lindstrom} and Brenti \cite[Theorem 4.2]{brenti} show that a real matrix $A$ is totally nonnegative if and only if $A$ is the weighted path matrix of a weighted planar network.  A good introduction to these ideas is given in 
\cite{fz}. 

In this section, we illustrate the results of the earlier sections by considering Holte's ``Amazing Matrix'' and show that various cases are indeed totally nonnegative, as conjectured in \cite{df}.

Holte's ``Amazing Matrix'' is the $n\times n$ matrix $P=(p_{ij})$ with 
\[
p_{ij}=\frac{1}{b^n}\sum_{r=0}^{j-\lfloor i/b\rfloor}(-1)^r \binom{n+1}{r}\binom{n-1-i+(j+1-r)b}{n}.
\]
The entry $p_{ij}$ gives the probability that, when adding $n$ random numbers in base $b$, the next carry will be $j$, given that the previous carry was $i$, see \cite{holte}. (Note that Holte indexes the rows and columns of an $n\times n$ amazing matrix by $0,1,\dots,n-1$.) Holte shows that $P$ is cross-symmetric \cite[Theorem 2]{holte}. It will be convenient for us to use the scaled version $P':=b^nP$ to avoid fractions. As examples of $P'$, when $b=3$ the $2\times 2, 3\times 3$ and $4\times 4$ matrices are displayed below.
\[
\begin{pmatrix}
6&3\\3&6
\end{pmatrix}
\qquad
\begin{pmatrix}
10&16&1\\
4&19&4\\
1&16&10
\end{pmatrix}
\qquad
 \begin{pmatrix}
15&51&15&0\\
5&45&30&1\\
1&30&45&5\\
0&15&51&15
\end{pmatrix}
\]

The factorisation of the $2\times 2$ matrix above that is given by the algorithm is:
\[
\begin{pmatrix}
6&3\\3&6
\end{pmatrix}
=
\begin{pmatrix}
4/3&2/3\\
2/3&4/3
\end{pmatrix}
 \begin{pmatrix}
9/2&0\\0&9/2
\end{pmatrix}
\]
and the first factor is an atom of the second kind where $c=1/2$.

The factorisation into atoms of the $n=b=3$ scaled amazing matrix  that is given by the algorithm developed earlier is:

\[
\begin{pmatrix}
10&16&1\\
4&19&4\\
1&16&10
\end{pmatrix}
=
\begin{pmatrix}
1&1/4&0\\
0&1&0\\
0&1/4&1
\end{pmatrix}
\begin{pmatrix}
1&0&0\\
4/9&1&4/9\\
0&0&1
\end{pmatrix}
\begin{pmatrix}
1&5/4&0\\
0&1&0\\
0&5/4&1
\end{pmatrix}
\begin{pmatrix}
9&0&0\\
0&9&0\\
0&0&9
\end{pmatrix},
\]
verifying that the matrix is totally nonnegative, as each of the factors is totally nonnegative.

The above factorisation gives rise to the  planar network in Figure~\ref{figure-3x3} whose weighted path matrix is the amazing matrix with $n=b=3$ (edges without weights attached are taken to have weight $1$, and all edges are directed from left to right): 
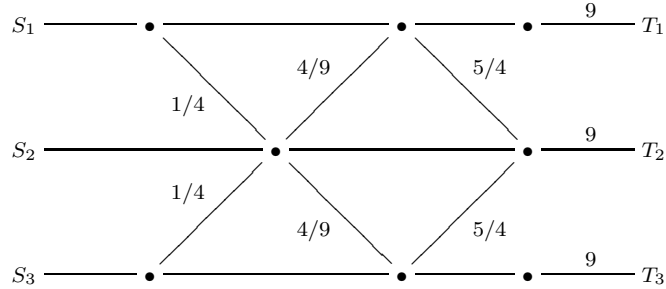
\begin{figure}[H]
\centering
$
\xymatrixrowsep{10ex}\xymatrixcolsep{10ex}\def\objectstyle{\scriptstyle}
\xymatrix@!0{
 S_1\edge[r]&\bullet\edge[rr]\edge[dr]_{1/4}&\edge[r]&\bullet\edge[r]\edge[dr]^{5/4}&\bullet\edge[r] ^{9}&T_1\\
 S_2\edge[rr]&\edge[r]&\bullet\edge[rr]\edge[ru]^{4/9}\edge[rd]_{4/9}&\edge[r]&\bullet\edge[r] ^{9}&T_2\\
  S_3\edge[r]&\bullet\edge[rr]\edge[ru]^{1/4}&\edge[r]&\bullet\edge[r]\edge[ur]_{5/4}&\bullet\edge[r] ^{9}&T_3\\
}$
\caption{\label{figure-3x3}Weighted planar networks for  scaling of the $3\times 3$ Amazing Matrix with base $b=3$ 
(all edges directed left to right)
 }
\end{figure}

The cross-symmetry property is manifested in the planar network by the symmetry in the horizontal line at the centre of the network. \\

The factorisation of the scaled amazing matrix in the case where $n=4$ and $b=3$ contains atoms of the second kind and leads to the planar network shown in Figure~\ref{figure-n=4,b=3}. Again, total nonnegativity is verified and the cross-symmetry property is seen by the symmetry about the horizontal line through the middle of the network.

\begin{figure}[H]
\centering
$
\xymatrixrowsep{10ex}\xymatrixcolsep{10ex}\def\objectstyle{\scriptstyle}
\xymatrix@!0{
 S_1\edge[rrrr]&&&&\bullet\edge[rr]\edge[dr]^{\frac{25}{27}}&&\bullet\edge[rr]^{15}\edge[dr]_{\frac{875}{324}}&&T_1\\
 S_2\edge[rr]&\bullet\edge[rr]^{\frac{5}{6}}\edge[dr]_{\frac{5}{24}}&&\bullet\edge[rr]\edge[ur]^{\frac{8}{25}}&&\bullet\edge[rr]_{\frac{7}{12}}\edge[dr]_{\frac{35}{24}}&&\bullet\edge[r]_{\frac{1944}{175}}&T_2\\
 &&\bullet\edge[ur]\edge[dr]&&&&\bullet\edge[ur]\edge[dr]&&\\
 S_3\edge[rr]&\bullet\edge[rr]_{\frac{5}{6}}\edge[ur]^{\frac{5}{24}}&&\bullet\edge[rr]\edge[dr]_{\frac{8}{25}}&&\bullet\edge[rr]^{\frac{7}{12}}\edge[ur]^{\frac{35}{24}}&&\bullet\edge[r]^{\frac{1944}{175}}&T_3\\
 S_4\edge[rrrr]&&&&\bullet\edge[rr]\edge[ur]_{\frac{25}{27}}&&\bullet\edge[rr]^{15}\edge[ur]^{\frac{875}{324}}&&T_4
}$
\caption{\label{figure-n=4,b=3}Weighted planar network for  scaling of the $4\times 4$ Amazing Matrix with base $b=3$
(all edges directed left to right) 
 }
\end{figure}
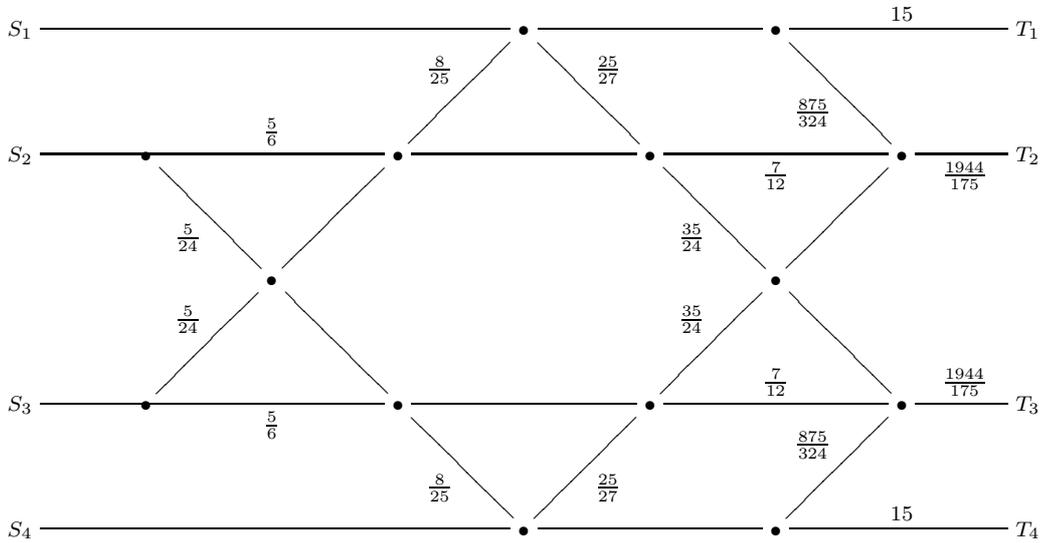

\begin{remark}
We have run the algorithm in arbitrary base $b$ up to $n=6$ and the total nonnegativity conjecture is verified in each of these cases. 
\end{remark}





\vskip 1cm

\noindent Maxwell Institute,\\
School of Mathematics,\\ University of Edinburgh,\\
James Clerk Maxwell Building,\\
The King's Buildings,\\
Peter Guthrie Tait Road,\\
Edinburgh EH9 3FD \\
           UK\\[0.5ex]
           
 \noindent           tom@maths.ed.ac.uk\\
 \noindent           andrewneate13@gmail.com

\end{document}